\documentclass{amsart}
\usepackage{amsmath}
\usepackage{amstext}
\usepackage{amssymb}
\usepackage{amsthm}
\swapnumbers
\theoremstyle{plain}
\newtheorem{thm}{Theorem}[section]
\newtheorem{lem}[thm]{Lemma}
\newtheorem{prop}[thm]{Proposition}
\newtheorem{cor}[thm]{Corollary}

\theoremstyle{definition}
\newtheorem{defi}[thm]{Definition}
\newtheorem{defs}[thm]{Definitions}
\newtheorem{defrmks}[thm]{Definition and Remarks}

\newtheorem{ntn}[thm]{Notation}
\newtheorem{rmd}[thm]{Reminder}

\theoremstyle{remark}
\newtheorem*{note}{Note}
\newtheorem{rmk}[thm]{Remark}

 \DeclareMathOperator{\height}{ht}
 
\DeclareMathOperator{\Ext}{Ext} \DeclareMathOperator{\Hom}{Hom}
 
\DeclareMathOperator{\Ker}{Ker}

\DeclareMathOperator{\Spec}{Spec} 
 \DeclareMathOperator{\ann}{ann}
\DeclareMathOperator{\grann}{gr-ann}

\def\Z{\mathbb Z}
\def\N{\mathbb N}

\def\fa{{\mathfrak{a}}}

\def\fb{{\mathfrak{b}}}
\def\fB{{\mathfrak{B}}}
\def\fc{{\mathfrak{c}}}
\def\fd{{\mathfrak{d}}}

\def\fm{{\mathfrak{m}}}

\def\fp{{\mathfrak{p}}}

\def\fq{{\mathfrak{q}}}

\def\fr{{\mathfrak{r}}}

\def\nn{\relax\ifmmode{\mathbb N_{0}}\else$\mathbb N_{0}$\fi}
\def\lra{\longrightarrow}

\begin{document}

\title[GRADED ANNIHILATORS AND TIGHT CLOSURE]{GRADED ANNIHILATORS OF MODULES OVER
THE FROBENIUS SKEW POLYNOMIAL RING, AND TIGHT CLOSURE}
\author{RODNEY Y. SHARP}
\address{Department of Pure Mathematics,
University of Sheffield, Hicks Building, Sheffield S3 7RH, United Kingdom\\
{\it Fax number}: 0044-114-222-3769}
\email{R.Y.Sharp@sheffield.ac.uk}

\thanks{The author was partially supported by the
Engineering and Physical Sciences Research Council of the United
Kingdom (Overseas Travel Grant Number EP/C538803/1).}


\date{\today}


\begin{abstract}
This paper is concerned with the tight closure of an ideal $\fa$
in a commutative Noetherian local ring $R$ of prime characteristic
$p$. Several authors, including R. Fedder, K.-i. Watanabe, K. E.
Smith, N. Hara and F. Enescu, have used the natural Frobenius
action on the top local cohomology module of such an $R$ to good
effect in the study of tight closure, and this paper uses that
device. The main part of the paper develops a theory of what are
here called `special annihilator submodules' of a left module over
the Frobenius skew polynomial ring associated to $R$; this theory
is then applied in the later sections of the paper to the top
local cohomology module of $R$ and used to show that, if $R$ is
Cohen--Macaulay, then it must have a weak parameter test element,
even if it is not excellent.
\end{abstract}

\maketitle

\setcounter{section}{-1}
\section{\sc Introduction}
\label{in}

Throughout the paper, $R$ will denote a commutative
Noetherian ring of prime characteristic $p$.
We shall always
denote by $f:R\lra R$ the Frobenius homomorphism, for which $f(r)
= r^p$ for all $r \in R$. Let $\fa$ be an ideal of $R$. The {\em
$n$-th Frobenius power\/} $\fa^{[p^n]}$ of $\fa$ is the ideal of
$R$ generated by all $p^n$-th powers of elements of $\fa$.

We use $R^{\circ}$ to denote the complement in $R$ of the union of
the minimal prime ideals of $R$. An element $r \in R$ belongs to
the {\em tight closure $\fa^*$ of $\fa$\/} if and only if there
exists $c \in R^{\circ}$ such that $cr^{p^n} \in \fa^{[p^n]}$ for
all $n \gg 0$. We say that $\fa$ is {\em tightly closed\/}
precisely when $\fa^* = \fa$. The theory of tight closure was
invented by M. Hochster and C. Huneke \cite{HocHun90}, and many
applications have been found for the theory: see \cite{Hunek96}
and \cite{Hunek98}, for example.

In the case when $R$ is local, several authors have used, as an
aid to the study of tight closure, the natural Frobenius action on
the top local cohomology module of $R$: see, for example, R.
Fedder \cite{F87}, Fedder and K.-i. Watanabe \cite{FW87}, K. E.
Smith \cite{S94}, N. Hara and Watanabe \cite{HW96} and F. Enescu
\cite{Enesc03}. This device is employed in this paper. The natural
Frobenius action provides the top local cohomology module of $R$
with a natural structure as a left module over the skew polynomial
ring $R[x,f]$ associated to $R$ and $f$. Sections \ref{nt} and
\ref{ga} develop a theory of what are here called `special
annihilator submodules' of a left $R[x,f]$-module $H$. To explain
this concept, we need the definition of the {\em graded
annihilator\/} $\grann_{R[x,f]}H$ of $H$. Now $R[x,f]$ has a
natural structure as a graded ring, and $\grann_{R[x,f]}H$ is
defined to be the largest graded two-sided ideal of $R[x,f]$ that
annihilates $H$. On the other hand, for a graded two-sided ideal
$\fB$ of $R[x,f]$, the {\em annihilator of $\fB$ in $H$\/} is
defined as
$$ \ann_H\fB := \{ h \in H : \theta h = 0 \mbox{~for all~}\theta
\in \fB\}. $$ I say that an $R[x,f]$-submodule of $H$ is a {\em
special annihilator submodule\/} of $H$ if it has the form
$\ann_H\fB$ for some graded two-sided ideal $\fB$ of $R[x,f]$.

There is a natural bijective inclusion-reversing correspondence
between the set of all special annihilator submodules of $H$ and
the set of all graded annihilators of submodules of $H$. A large
part of this paper is concerned with exploration and exploitation
of this correspondence. It is particularly satisfactory in the
case where the left $R[x,f]$-module $H$ is $x$-torsion-free, for
then it turns out that the set of all graded annihilators of
submodules of $H$ is in bijective correspondence with a certain
set of radical ideals of $R$, and one of the main results of \S
\ref{ga} is that this set is finite in the case where $H$ is
Artinian as an $R$-module. The theory that emerges has some
uncanny similarities to tight closure theory. Use is made of the
Hartshorne--Speiser--Lyubeznik Theorem (see R. Hartshorne and R.
Speiser \cite[Proposition 1.11]{HarSpe77}, G. Lyubeznik
\cite[Proposition 4.4]{Lyube97}, and M. Katzman and R. Y. Sharp
\cite[1.4 and 1.5]{KS}) to pass between a general left
$R[x,f]$-module that is Artinian over $R$ and one that is
$x$-torsion-free.

In \S \ref{tc}, this theory of special annihilator submodules is
applied to prove an existence theorem for weak parameter test elements
in a Cohen--Macaulay local ring of characteristic $p$. To explain this,
I now review some definitions concerning weak test elements.

A {\em $p^{w_0}$-weak test element\/} for $R$ (where $w_0$ is a
non-negative integer) is an element $c' \in R^{\circ}$ such that,
for every ideal $\fb$ of $R$ and for $r \in R$, it is the case
that $r \in \fb^*$ if and only if $c'r^{p^n} \in \fb^{[p^n]}$ for
all $n \geq w_0$. A $p^0$-weak test element is called a {\em test
element\/}.

A proper ideal $\fa$ in $R$ is said to be a {\em parameter ideal\/}
precisely when it can be generated by $\height \fa$ elements. Parameter ideals
play an important r\^ole in tight closure theory, and Hochster and Huneke
introduced the concept of parameter test element for $R$.
A {\em $p^{w_0}$-weak parameter test element\/} for $R$
is an element $c' \in R^{\circ}$ such that,
for every parameter ideal $\fb$ of $R$ and for $r \in R$, it is the case
that $r \in \fb^*$ if and only if $c'r^{p^n} \in \fb^{[p^n]}$ for
all $n \geq w_0$. A $p^0$-weak parameter test element is called a {\em
parameter test
element\/}.

It is a result of Hochster and Huneke \cite[Theorem
(6.1)(b)]{HocHun94} that an algebra of finite type over an
excellent local ring of characteristic $p$ has a $p^{w_0}$-weak
test element for some non-negative integer $w_0$; furthermore,
such an algebra which is reduced actually has a test element. Of
course, a (weak) test element is a (weak) parameter test element.

One of the main results of this paper is Theorem \ref{tc.2}, which
shows that every Cohen--Macaulay local ring of characteristic $p$,
even if it is not excellent, has a $p^{w_0}$-weak parameter test
element for some non-negative integer $w_0$.

Lastly, the final \S \ref{en} establishes some connections between
the theory developed in this paper and the $F$-stable primes of F. Enescu
\cite{Enesc03}.

\section{\sc Graded annihilators and related concepts}
\label{nt}

\begin{ntn}
\label{nt.1} Throughout, $R$ will denote a commutative Noetherian
ring of prime characteristic $p$. We shall work with the
 skew polynomial ring $R[x,f]$ associated to $R$ and $f$ in the
indeterminate $x$ over $R$. Recall that $R[x,f]$ is, as a left
$R$-module, freely generated by $(x^i)_{i \in \nn}$ (I use $\N$
and $\nn$ to denote the set of positive integers and the set of
non-negative integers, respectively),
 and so consists
 of all polynomials $\sum_{i = 0}^n r_i x^i$, where  $n \in \nn$
 and  $r_0,\ldots,r_n \in R$; however, its multiplication is subject to the
 rule
 $$
  xr = f(r)x = r^px \quad \mbox{~for all~} r \in R\/.
 $$
Note that $R[x,f]$ can be considered as a positively-graded ring
$R[x,f] = \bigoplus_{n=0}^{\infty} R[x,f]_n$, with $R[x,f]_n =
Rx^n$ for all $n \in \nn$. The ring $R[x,f]$ will be referred to as the {\em
Frobenius skew polynomial ring over $R$.}

Throughout, we shall let $G$ and $H$ denote left $R[x,f]$-modules.
The {\em annihilator of $H$\/} will be denoted by $\ann_{R[x,f]}H$
or $\ann_{R[x,f]}(H)$. Thus
$$
\ann_{R[x,f]}(H) = \{ \theta \in R[x,f] : \theta h = 0 \mbox{~for all~} h \in H\},
$$
and this is a (two-sided) ideal of $R[x,f]$. For a two-sided ideal
$\fB$ of $R[x,f]$, we shall use $\ann_H\fB$ or $\ann_H(\fB)$ to
denote the {\em annihilator of $\fB$ in $H$}. Thus
$$
\ann_H\fB = \ann_H(\fB) = \{ h \in
H : \theta h = 0 \mbox{~for all~}\theta \in \fB\},
$$
and this is an $R[x,f]$-submodule of $H$.
\end{ntn}

\begin{defrmks}
\label{nt.2}
We say that the left $R[x,f]$-module $H$ is {\em $x$-torsion-free\/} if $xh = 0$, for
$h \in H$, only when $h = 0$.
The set $\Gamma_x(H) := \left\{ h \in H : x^jh = 0
\mbox{~for some~} j \in \N \right\}$ is an $R[x,f]$-submodule of
$H$, called the\/ {\em $x$-torsion submodule} of $H$. The
$R[x,f]$-module $H/\Gamma_x(H)$ is $x$-torsion-free.
\end{defrmks}

\begin{rmk}
\label{nt.2b} Let $\fB$ be a subset of $R[x,f]$. It is easy to
see that $\fB$ is a graded two-sided ideal of $R[x,f]$ if and only if
there is an ascending chain $(\fb_n)_{n \in \nn}$ of ideals of $R$ (which must,
of course, be eventually stationary) such that
$\fB = \bigoplus_{n\in\nn}\fb_n x^n$. We shall sometimes denote the
ultimate constant value of the ascending sequence $(\fb_n)_{n \in \nn}$ by
$\lim_{n \rightarrow \infty}\fb_n$.

Note that, in particular, if $\fb$ is an
ideal of $R$, then $\fb R[x,f] = \bigoplus_{n \in \nn} \fb x^n$ is
a graded two-sided ideal of $R[x,f]$. It was noted in \ref{nt.1}
that the annihilator of a left $R[x,f]$-module is a two-sided
ideal.
\end{rmk}

\begin{lem} [Y. Yoshino {\cite[Corollary (2.7)]{Yoshi94}}]
\label{nt.2c} The ring $R[x,f]$ satisfies the ascending chain condition on
graded two-sided ideals.
\end{lem}

\begin{proof}  This can be proved by the argument in Yoshino's proof of
\cite[Corollary (2.7)]{Yoshi94}.
\end{proof}

\begin{defs}
\label{nt.2d} We define the {\em graded annihilator\/} $\grann_{R[x,f]}H$ of
the left $R[x,f]$-module $H$ by
$$
\grann_{R[x,f]}H = \left\{ \sum_{i=0}^n r_ix^i \in R[x,f] : n \in
\nn, \mbox{~and~} r_i \in R,\, r_ix^i \in \ann_{R[x,f]}H \mbox{~for
all~} i = 0, \ldots, n\right\}.
$$
Thus $\grann_{R[x,f]}H$ is the largest graded two-sided ideal of $R[x,f]$
contained in $\ann_{R[x,f]}H$; also, if we write $\grann_{R[x,f]}H =
\bigoplus_{n\in\nn}\fb_n x^n$ for a suitable
ascending chain $(\fb_n)_{n \in \nn}$ of ideals of $R$, then $\fb_0 = (0:_RH)$,
the annihilator of $H$ as an $R$-module.

We say that an $R[x,f]$-submodule of $H$ is a {\em special
annihilator submodule of $H$\/} if it has the form $\ann_H(\fB)$
for some {\em graded\/} two-sided ideal $\fB$ of $R[x,f]$. We
shall use $\mathcal{A}(H)$ to denote the set of special
annihilator submodules of $H$.
\end{defs}

\begin{defrmks}
\label{nt.05d} There are some circumstances in which $\grann_{R[x,f]}H =
\ann_{R[x,f]}H$: for example, this would be the case if $H$ was a $\Z$-graded
left $R[x,f]$-module. Work of Y. Yoshino in \cite[\S 2]{Yoshi94} provides us
with further examples.

Following Yoshino \cite[Definition (2.1)]{Yoshi94},
we say that $R$ {\em has sufficiently many units\/} precisely when,
for each $n \in \N$, there exists $r_n \in R$ such that all $n$ elements
$(r_n)^{p^i} - r_n~(1 \leq i \leq n)$ are units of $R$.
Yoshino proved in \cite[Lemma (2.2)]{Yoshi94} that if
either $R$ contains an infinite field,
or $R$ is local and has infinite residue field,
then $R$ has sufficiently many units. He went on to show in
\cite[Theorem (2.6)]{Yoshi94} that, if
$R$ has sufficiently many units,
then each two-sided ideal
of $R[x,f]$ is graded.

Thus if $R$ has sufficiently many units, then $\grann_{R[x,f]}H =
\ann_{R[x,f]}H$, even if $H$ is not graded.
\end{defrmks}

\begin{lem}
\label{nt.3} Let $\fB$ and $\fB '$ be
graded two-sided ideals of $R[x,f]$ and let $N$ and
$N'$ be $R[x,f]$-submodules of the left $R[x,f]$-module $H$.

\begin{enumerate}
\item If $\fB \subseteq \fB '$, then $\ann_H(\fB) \supseteq \ann_H(\fB ')$.
\item If $N \subseteq N'$, then $\grann_{R[x,f]}N \supseteq \grann_{R[x,f]}N '$.
\item We have
$\fB \subseteq \grann_{R[x,f]}\left(\ann_H(\fB)\right)$.
\item We have $N \subseteq \ann_H\!\left(\grann_{R[x,f]}N\right).$
\item There
is an order-reversing bijection, $\Gamma\/,$ from the set
$\mathcal{A}(H)$ of special annihilator submo\-d\-ules of $H$ to
the set of graded annihilators of submodules of $H$ given by
$$
\Gamma : N \longmapsto \grann_{R[x,f]}N.
$$
The inverse bijection, $\Gamma^{-1},$ also order-reversing, is
given by
 $$
\Gamma^{-1} : \fB \longmapsto \ann_H(\fB).
$$
\end{enumerate} \end{lem}

\begin{proof} Parts (i), (ii), (iii) and (iv) are obvious.

(v) Application of part (i) to the inclusion in part (iii) yields that
$$
\ann_H(\fB) \supseteq \ann_H\!\left(
\grann_{R[x,f]}\left(\ann_H(\fB)\right) \right)\mbox{;}
$$
however, part (iv) applied to the $R[x,f]$-submodule
$\ann_H(\fB)$ of $H$ yields that
$$
\ann_H(\fB) \subseteq \ann_H\!\left(
\grann_{R[x,f]}\left(\ann_H(\fB)\right) \right)\mbox{;}
$$
hence $ \ann_H(\fB) = \ann_H\!\left(
\grann_{R[x,f]}\left(\ann_H(\fB)\right) \right). $ Similar
considerations show that
$$
\grann_{R[x,f]}N =
\grann_{R[x,f]}\left(\ann_H\!\left(\grann_{R[x,f]}N\right)\right).
$$
\end{proof}

\begin{rmk}
\label{nt.4} It follows from Lemma \ref{nt.3} that, if $N$ is a special annihilator
submodule of $H$, then it is the annihilator (in $H$) of its own
graded annihilator. Likewise,
a graded two-sided ideal $\fB$ of $R[x,f]$ which is the
graded annihilator of some $R[x,f]$-submodule
of $H$ must be the graded annihilator of $\ann_H(\fB)$.
\end{rmk}

Much use will be made of the
following lemma.

\begin{lem}
\label{nt.5} Assume that the left $R[x,f]$-module $G$ is
$x$-torsion-free.

Then there is a radical ideal $\fb$ of $R$ such that
$\grann_{R[x,f]}G = \fb R[x,f] = \bigoplus _{n\in\nn} \fb x^n$.
\end{lem}

\begin{proof} There is a family
$(\fb_n)_{n\in\nn}$ of ideals of $R$ such that $\fb_n \subseteq
\fb_{n+1}$ for all $n \in \nn$ and $\grann_{R[x,f]}G = \bigoplus
_{n\in\nn} \fb_n x^n$. There exists $n_0 \in \nn$ such that $\fb_n = \fb_{n_0}$
for all $n \geq n_0$. Set $\fb := \fb_{n_0}$. It is enough for us to show that,
if $r \in R$ and $e \in \nn$ are such that $r^{p^e} \in \fb$, then $r \in \fb_0$.

To this end, let $h \in \N$ be such that $h \geq \max \{e,n_0\}$. Then, for all
$g \in G$, we have
$x^hrg = r^{p^h}x^hg = 0$, since $r^{p^h} \in \fb = \fb_h$. Since $G$ is
$x$-torsion-free, it follows that $rG = 0$, so that $r \in \fb_0$.
\end{proof}

\begin{defi}
\label{nt.6} Assume that the left $R[x,f]$-module $G$ is
$x$-torsion-free.  An ideal $\fb$ of $R$ is called a {\em
$G$-special $R$-ideal} if there is an $R[x,f]$-submodule
$N$ of $G$ such that $\grann_{R[x,f]}N = \fb R[x,f] = \bigoplus
_{n\in\nn} \fb x^n$. It is worth noting that, then, the ideal
$\fb$ is just $(0:_RN)$.

We shall denote the set of $G$-special $R$-ideals by
$\mathcal{I}(G)$. Note that, by Lemma \ref{nt.5}, all the ideals
in $\mathcal{I}(G)$ are radical.
\end{defi}

We can now combine together the results of Lemmas \ref{nt.3}(v)
and \ref{nt.5} to obtain the following result, which is
fundamental for the work in this paper.

\begin{prop}
\label{nt.7}  Assume that
the left $R[x,f]$-module $G$ is $x$-torsion-free.

There is an order-reversing bijection, $\Delta : \mathcal{A}(G)
\lra \mathcal{I}(G),$ from the set $\mathcal{A}(G)$ of special
annihilator submodules of $G$ to the set $\mathcal{I}(G)$ of
$G$-special $R$-ideals given by
$$
\Delta : N \longmapsto \left(\grann_{R[x,f]}N\right)\cap R =
(0:_RN).
$$

The inverse bijection, $\Delta^{-1} : \mathcal{I}(G) \lra
\mathcal{A}(G),$ also order-reversing, is given by
$$
\Delta^{-1} : \fb \longmapsto \ann_G\left(\fb R[x,f])\right).
$$
When $N \in \mathcal{A}(G)$ and $\fb \in \mathcal{I}(G)$ are such
that $\Delta(N) = \fb$, we shall say simply that `{\em $N$ and
$\fb$ correspond}'.
\end{prop}

\begin{cor}
\label{nt.9} Assume that
the left $R[x,f]$-module $G$ is $x$-torsion-free.

Then both the sets $\mathcal{A}(G)$ and $\mathcal{I}(G)$ are
closed under taking arbitrary intersections.
\end{cor}

\begin{proof} Let $\left(N_{\lambda}\right)_{\lambda \in \Lambda}$ be an
arbitrary family of special annihilator submodules of $G$. For
each $\lambda \in \Lambda$, let $\fb_{\lambda}$ be the
$G$-special $R$-ideal corresponding to $N_{\lambda}$. In
view of Proposition \ref{nt.7}, it is sufficient for us to show
that $\bigcap_{\lambda \in \Lambda}N_{\lambda} \in \mathcal{A}(G)$
and $\fb := \bigcap_{\lambda \in \Lambda}\fb_{\lambda} \in
\mathcal{I}(G)$.

To prove these, simply note that
$$
{\textstyle \bigcap_{\lambda \in \Lambda}N_{\lambda} =
\bigcap_{\lambda \in \Lambda}\ann_G\left(\fb_{\lambda}
R[x,f])\right) = \ann_G\left(\left(\sum_{\lambda \in
\Lambda}\fb_{\lambda}\right)\! R[x,f]\right)}
$$
and that $\sum_{\lambda \in \Lambda}N_{\lambda}$ is an
$R[x,f]$-submodule of $G$ such that
$$
{\textstyle \grann_{R[x,f]}\left(\sum_{\lambda \in
\Lambda}N_{\lambda}\right) = \bigcap_{\lambda \in
\Lambda}\grann_{R[x,f]}N_{\lambda} = \bigcap_{\lambda \in
\Lambda}\left(\fb_{\lambda}R[x,f] \right) = \fb R[x,f].}
$$
\end{proof}

\begin{rmk}
\label{nt.8}  Suppose
that the left $R[x,f]$-module $G$ is $x$-torsion-free.

It is worth pointing out now that, since $R$ is Noetherian, so
that the set $\mathcal{I}(G)$ of $G$-special $R$-ideals
satisfies the ascending chain condition, it is a consequence of
Proposition \ref{nt.7} that the set $\mathcal{A}(G)$ of special
annihilator submodules of $G$, partially ordered by inclusion,
satisfies the descending chain condition. This is the case even if
$G$ is not finitely generated.  Note that (by \cite[Theorem
(1.3)]{Yoshi94}), the (noncommutative) ring $R[x,f]$ is neither
left nor right Noetherian if $\dim R > 0$.
\end{rmk}

\section{\sc Examples relevant to the theory of tight closure}
\label{ex}

The purpose of this section is to present some motivating examples, from the
theory of tight closure, of some of the concepts introduced in \S \ref{nt}.
Throughout this section, we shall again employ the notation of
\ref{nt.1}, and $\fa$ will always denote an ideal of $R$.
Recall that the {\em Frobenius closure $\fa^F$ of $\fa$} is the ideal of $R$
defined by
$$ \fa ^F := \big\{ r
\in R : \mbox{there exists~} n \in \nn \mbox{~such that~} r^{p^n}
\in \fa^{[p^n]}\big\}.
$$

\begin{rmk}
\label{ex.0} Let $(\fb_n)_{n\in\nn}$ be a family of ideals of $R$
such that $\fb_n \subseteq f^{-1}(\fb_{n+1})$ for all $n \in \nn$.
Then $\bigoplus_{n\in\nn}\fb_nx^n$ is a graded left ideal of
$R[x,f]$, and so we may form the graded left $R[x,f]$-module $
R[x,f]/\bigoplus_{n\in\nn}\fb_nx^n$. This may be viewed as
$\bigoplus_{n\in\nn}R/\fb_n$, where, for $r \in R$ and $n \in
\nn$, the result of multiplying the element $r + \fb_n$ of the
$n$-th component by $x$ is the element $r^{p} + \fb_{n+1}$ of the
$(n+1)$-th component.

Note that the left $R[x,f]$-module
$R[x,f]/\bigoplus_{n\in\nn}\fb_nx^n$ is $x$-torsion-free if and
only if $\fb_n = f^{-1}(\fb_{n+1})$ for all $n \in \nn$, that is,
if and only if $(\fb_n)_{n\in\nn}$ is an $f$-sequence in the sense of
\cite[Definition 4.1(ii)]{SN}.
\end{rmk}

\begin{ntn}
\label{ex.0a} Since $R[x,f]\fa =
\bigoplus_{n\in\nn}\fa^{[p^n]}x^n$, we can view the graded left
$R[x,f]$-module $$R[x,f]/R[x,f]\fa$$ as
$\bigoplus_{n\in\nn}R/\fa^{[p^n]}$ in the manner described in
\ref{ex.0}. We shall denote the graded left $R[x,f]$-module
$\bigoplus_{n\in\nn}R/\fa^{[p^n]}$ by $H(\fa)$.

Recall from \cite[4.1(iii)]{SN} that
$\left((\fa^{[p^n]})^F\right)_{n \in \nn}$ is the {\em canonical
$f$-sequence associated to $\fa$}. We shall denote
$\bigoplus_{n\in\nn}R/(\fa^{[p^n]})^F$, considered as a graded
left $R[x,f]$-module in the manner described in \ref{ex.0}, by
$G(\fa)$. Note that $G(\fa)$ is $x$-torsion-free.
\end{ntn}

\begin{lem}
\label{ex.0b} With the notation of\/ {\rm \ref{ex.0a}}, we have
$\Gamma_x(H(\fa)) = \bigoplus_{n\in\nn}(\fa^{[p^n]})^F/\fa^{[p^n]}$,
so that there is an isomorphism of graded left $R[x,f]$-modules
$$
H(\fa)/ \Gamma_x(H(\fa)) \cong G(\fa).
$$
\end{lem}

\begin{proof} Let $n \in \nn$ and $r \in R$.
Then the element $r + \fa^{[p^n]}$ of the $n$-th
component of $H(\fa)$ belongs to
$\Gamma_x(H(\fa))$
if and only if there exists $m \in \nn$ such
that $x^m(r + \fa^{[p^n]})= r^{p^m} + (\fa^{[p^{n}]})^{[p^{m}]} = 0$, that is,
if and only if $r \in (\fa^{[p^n]})^F$.
\end{proof}

\begin{prop}
\label{ex.1} We use the notation of\/ {\rm \ref{ex.0a}}.

Suppose that there exists a $p^{w_0}$-weak test element $c$ for
$R$, for some $w_0 \in \nn$. Then
\begin{enumerate}
\item $ \ann_{H(\fa)}\left(
\bigoplus_{n\geq w_0} Rcx^n\right) =
\bigoplus_{n\in\nn}(\fa^{[p^n]})^*/\fa^{[p^n]}\mbox{;} $

\item $
\ann_{G(\fa)}\left(\bigoplus_{n\in
\nn} Rcx^n\right) =
\bigoplus_{n\in\nn}(\fa^{[p^n]})^*/(\fa^{[p^n]})^F. $
\end{enumerate}
\end{prop}

\begin{proof} (i)  Let $j \in \nn$ and $r \in R$.
Then the element $r + \fa^{[p^j]}$ of the $j$-th
component of $H(\fa)$ belongs to
$\ann_{H(\fa)}\left(\bigoplus_{n\geq w_0} Rcx^n\right)$ if and only if $cr^{p^n} \in
(\fa^{[p^j]})^{[p^n]}$ for all $n \geq w_0$, that is, if and only if $r \in
(\fa^{[p^j]})^*$.

(ii) By part (i),
$$
{\textstyle \bigoplus_{n\in \nn}}(\fa^{[p^n]})^*/(\fa^{[p^n]})^F
\subseteq \ann_{G(\fa)}\left({\textstyle \bigoplus_{n\geq w_0}}
Rcx^n\right).
$$
Note that $\ann_{G(\fa)}\left(\bigoplus_{n\geq w_0} Rcx^n\right)$
is a graded $R[x,f]$-submodule
of $G(\fa)$. Let $j \in \nn$ and $r \in R$ be such that
$r + (\fa^{[p^j]})^F$ belongs to the $j$-th
component of $\ann_{G(\fa)}\left(\bigoplus_{n\geq w_0} Rcx^n\right)$. Then,
for all $n \geq w_0$, we have $cr^{p^n} \in (\fa^{[p^{j+n}]})^F =
((\fa^{[p^{j}]})^{[p^{n}]})^F$. Therefore, by \cite[Lemma 0.1]{KS}, we
have $r \in (\fa^{[p^{j}]})^*$.

It follows from this that
$$
{\textstyle \bigoplus_{n\in \nn}}(\fa^{[p^n]})^*/(\fa^{[p^n]})^F
= \ann_{G(\fa)}\left({\textstyle \bigoplus_{n\geq w_0}} Rcx^n\right),
$$
and so $\bigoplus_{n\in \nn}(\fa^{[p^n]})^*/(\fa^{[p^n]})^F$ is a special annihilator
submodule of the $x$-torsion-free graded left $R[x,f]$-module $G(\fa)$.
Let $\fb$ be the $G(\fa)$-special $R$-ideal corresponding to
this member of $\mathcal{A}(G(\fa))$. The above-displayed equation shows that
$Rc \subseteq \fb$. Hence, by Proposition \ref{nt.7},
\begin{align*}
{\textstyle \bigoplus_{n\in \nn}}(\fa^{[p^n]})^*/(\fa^{[p^n]})^F &
= \ann_{G(\fa)}\left({\textstyle \bigoplus_{n\geq w_0}} Rcx^n\right)
\supseteq
\ann_{G(\fa)}\left({\textstyle \bigoplus_{n\in \nn}} Rcx^n\right) \\
& \supseteq \ann_{G(\fa)}\left({\textstyle \bigoplus_{n\in \nn}} \fb
x^n\right) = {\textstyle \bigoplus_{n\in
\nn}}(\fa^{[p^n]})^*/(\fa^{[p^n]})^F.
\end{align*}
\end{proof}

\begin{defi}
\label{ex.1t} The {\em weak test ideal $\tau'(R)$
of $R$\/} is defined to be
the ideal generated by $0$ and all weak test elements for $R$.  (By a
`weak test element' for $R$ we mean a $p^{w_0}$-weak
test element for $R$ for some $w_0 \in \nn$.)

It is easy to see that each element of $\tau'(R) \cap R^{\circ}$
is a weak test element for $R$.
\end{defi}

\begin{thm}
\label{ex.2} We use the notation of\/ {\rm \ref{ex.0a}}.
Suppose that there exists a $p^{w_0}$-weak test element $c$ for
$R$, for some $w_0 \in \nn$.

Let $H$ be the positively-graded left $R[x,f]$-module given by
$$
H := \bigoplus_{\fa \textit{~is an ideal of~} R}H(\fa) =
\bigoplus_{\fa \textit{~is an ideal of~}
R}\Big(\textstyle{\bigoplus_{n\in\nn}} R/\fa^{[p^n]}\Big).
$$
Set $T := {\displaystyle \bigoplus_{\fa \textit{~is an ideal of~}
R}} \Big(\bigoplus_{n\in\nn}(\fa^{[p^n]})^*/\fa^{[p^n]}\Big)$.

\begin{enumerate}
\item Then
$
T
= \ann_{H}\left(
\bigoplus_{n\geq w_0} Rcx^n\right),
$
and so is a special annihilator submodule of $H$.
\item Write $\grann_{R[x,f]}T =
\bigoplus_{n \in \nn} \fc_nx^n$ for a suitable ascending
chain $(\fc_n)_{n\in\nn}$ of ideals of $R$. Then $\lim_{n \rightarrow \infty}\fc_n
= \tau'(R)$, the weak test ideal for $R$.

\item Furthermore, $T$ contains every special annihilator
submodule $T'$ of $H$ for which the graded annihilator
$\grann_{R[x,f]}T' = \bigoplus_{n \in \nn} \fb_nx^n$ has\/
$\height (\lim_{n \rightarrow \infty}\fb_n) \geq 1$. (The height
of the improper ideal $R$ is considered to be $\infty$.)
\end{enumerate}
\end{thm}

\begin{proof} (i) This is immediate from Proposition \ref{ex.1}(i).

(ii) Write $\fc := \lim_{n \rightarrow \infty}\fc_n$.
Since there exists a weak test element for $R$, the ideal $\tau'(R)$
can be generated by finitely many weak test elements for $R$, say by
$c_i~(i = 1, \ldots, h)$, where $c_i$ is a $p^{w_i}$-weak test element for
$R$ (for $i= 1, \ldots, h$). Set $\widetilde{w} = \max\{w_1, \ldots, w_h\}$.
It is immediate from part (i) that $\bigoplus_{n \geq \widetilde{w}}\tau'(R)x^n
\subseteq \grann_{R[x,f]}T$, and so $\tau'(R) \subseteq \fc$. Therefore
$\height \fc \geq 1$, so that $\fc \cap R^{\circ} \neq \emptyset$ by
prime avoidance, and $\fc$ can be generated by its elements in $R^{\circ}$.

There exists $m_0 \in \nn$ such that $\fc_n = \fc$ for all $n \geq m_0$. Let
$c' \in \fc \cap R^{\circ}$. Thus $T$ is annihilated by $c'x^n$ for
all $n \geq m_0$; therefore, for each ideal $\fa$ of $R$, and for all $r \in \fa^*$,
we have $c'r^{p^n} \in \fa^{[p^n]}$ for all $n \geq m_0$, so that
$c'$ is a $p^{m_0}$-weak test element for $R$. Therefore $c' \in \tau'(R)$.
Since $\fc$ can be generated by elements in $\fc \cap R^{\circ}$, it
follows that $\fc \subseteq \tau'(R)$.

(iii) Since $T' = \ann_{H}\Big(\bigoplus_{n \in \nn} \fb_nx^n\Big)$, it follows that
$$
T' = \bigoplus_{\fa \text{~is an ideal of~} R} \Big({\textstyle
\bigoplus_{n\in\nn}}\fa_n/\fa^{[p^n]}\Big),
$$
where, for each ideal $\fa$ of $R$ and each $n \in \nn$, the ideal $\fa_n$ of
$R$ contains $\fa^{[p^n]}$. Suppose that $\lim_{n \rightarrow \infty}\fb_n = \fb$
and that $v_0 \in\nn$ is such
that $\fb_n = \fb$ for all $n \geq v_0$. Since $\height \fb \geq 1$, there
exists $\overline{c} \in \fb\cap R^{\circ}$, by prime avoidance.
Let $\fa$ be an ideal of $R$
and let $n \in \nn$.
Then, for each $r \in \fa_n$, the element $r + \fa^{[p^n]}$ of the
$n$-th component of $H(\fa)$ is annihilated by $\overline{c}x^j$ for all $j \geq
v_0$. This means that $\overline{c}r^{p^j} \in (\fa^{[p^n]})^{[p^j]}$ for all
$j \geq v_0$, so that $r \in (\fa^{[p^n]})^*$. Therefore $T' \subseteq T$.
\end{proof}

\begin{thm}
\label{ex.3} We use the notation of\/ {\rm \ref{ex.0a}}.
Suppose that there exists a $p^{w_0}$-weak test element $c$ for
$R$, for some $w_0 \in \nn$.

Let $G$ be the positively-graded $x$-torsion-free left $R[x,f]$-module given by
$$
G := \bigoplus_{\fa \textit{~is an ideal of~} R}G(\fa) =
\bigoplus_{\fa \textit{~is an ideal of~} R}\Big({\textstyle
\bigoplus_{n\in\nn}} R/(\fa^{[p^n]})^F\Big).
$$
Set $U := {\displaystyle \bigoplus_{\fa \textit{~is an ideal of~}
R}} \Big(\bigoplus_{n\in\nn}(\fa^{[p^n]})^*/(\fa^{[p^n]})^F\Big)$.

\begin{enumerate}
\item Then $ U = \ann_{G}\left( \bigoplus_{n\in\nn} Rcx^n\right),
$ and so is a special annihilator submodule of $G$. \item Let
$\fb$ be the $G$-special $R$-ideal corresponding to $U$. Then
$\fb$ is the smallest member of $\mathcal{I}(G)$ of positive
height.
\end{enumerate}
\end{thm}

\begin{proof} (i) This is immediate from Proposition \ref{ex.1}(ii).

(ii) Note that $Rc \subseteq \fb$, by part (i); therefore $\height
\fb \geq 1$. To complete the proof, we show that, if $\fb' \in
\mathcal{I}(G)$ has $\height \fb' \geq 1$, then $\fb \subseteq
\fb'$. By prime avoidance, there exists $\widetilde{c} \in \fb'
\cap R^{\circ}$. Let $U' \in \mathcal{A}(G)$ correspond to $\fb'$
(in the correspondence of Proposition \ref{nt.7}). Since $U' =
\ann_G\fb' R[x,f]$, it follows that
$$
U' = \bigoplus_{\fa \text{~is an ideal of~} R} \Big({\textstyle
\bigoplus_{n\in\nn}}\fa_n/(\fa^{[p^n]})^F\Big),
$$
where, for each ideal $\fa$ of $R$ and each $n \in \nn$, the ideal
$\fa_n$ of $R$ contains $(\fa^{[p^n]})^F$. Let $\fa$ be an ideal
of $R$ and let $n \in \nn$. Then, for each $r \in \fa_n$, the
element $r + (\fa^{[p^n]})^F$ of the $n$-th component of $G(\fa)$
is annihilated by $\widetilde{c}x^j$ for all $j \geq 0$. This
means that $\widetilde{c}r^{p^j} \in ((\fa^{[p^n]})^{[p^j]})^F$
for all $j \geq 0$, so that $r \in (\fa^{[p^n]})^*$ by \cite[Lemma
0.1(i)]{KS}. Therefore $U' \subseteq U$, so that $\fb' \supseteq
\fb$.
\end{proof}

\section{\sc Properties of
special annihilator submodules in the $x$-torsion-free case}
\label{ga}

Throughout this section, we shall employ the notation of
\ref{nt.1}. The aim is to develop the theory of special
annihilator submodules of an $x$-torsion-free left
$R[x,f]$-module.

\begin{lem}
\label{ga.1} Suppose that
$G$ is $x$-torsion-free. Let $N$ be a special annihilator submodule of
$G$.  Then the left $R[x,f]$-module $G/N$ is also
$x$-torsion-free.
\end{lem}

\begin{proof} By Lemma \ref{nt.5} and
Proposition \ref{nt.7}, there is a radical ideal $\fb$ of $R$ such
that $N = \ann_G\left(\fb R[x,f]\right).$ Let $g \in G$ be such
that $xg \in N$. Therefore, for all $r \in \fb$ and all $j \in
\nn$, we have $rx^j(xg) = 0$, that is $rx^{j+1}g = 0$. Also, for
$r \in \fb$, since $r(xg) = 0$, we have $x(rg) = r^pxg = 0$, and
so $rg = 0$ because $G$ is $x$-torsion-free. Thus $g \in
\ann_G\left(\bigoplus_{n\in\nn}\fb x^n\right) = N$. It follows
that $G/N$ is $x$-torsion-free.
\end{proof}

\begin{lem}
\label{ga.1a}  Suppose that $G$ is $x$-torsion-free. Let $\fa$ be
an ideal of $R$, and set $L := \ann_G\left(\fa R[x,f]\right) \in
\mathcal{A}(G)$. Then $L =
 \ann_G\left(\sqrt{\fa} R[x,f]\right)$.
\end{lem}

\begin{proof} Let $\fd \in \mathcal{I}(G)$ correspond to $L$.
Note that $\fd$ is radical, by Lemma \ref{nt.5}; also, $\fa
\subseteq \fd$. Hence $ \fa \subseteq \sqrt{\fa} \subseteq
\sqrt{\fd} = \fd$. Since $\ann_G\left(\fa R[x,f]\right) =
\ann_G\left(\fd R[x,f]\right)$, we must have $\ann_G\left(\fa
R[x,f]\right) = \ann_G\left(\sqrt{\fa} R[x,f]\right).$
\end{proof}

\begin{prop}
\label{ga.3} Suppose that $G$ is $x$-torsion-free. Let $\fa$ be an
ideal of $R$, and set $$L := \ann_G\left(\fa R[x,f]\right) \in
\mathcal{A}(G).$$ Note that $G/L$ is $x$-torsion-free, by Lemma\/
{\rm \ref{ga.1}}. Let $N$ be an $R$-submodule of $G$ such that $L
\subseteq N \subseteq G$.

\begin{enumerate}
\item If $N = \ann_G\left(\fb R[x,f]\right) \in \mathcal{A}(G)$,
where $\fb$ is an ideal of $R$ contained in $\fa$, then
$$N/L = \ann_{G/L}\left((\fb:\fa)R[x,f]\right) \in \mathcal{A}(G/L).$$
Furthermore, if the ideal in $\mathcal{I}(G)$ corresponding to $N$
is $\fb$, then $(\fb:\fa)$ is the ideal in $\mathcal{I}(G/L)$
corresponding to $N/L$.

\item If $N/L = \ann_{G/L}\left(\fc R[x,f]\right) \in
\mathcal{A}(G/L)$, where $\fc$ is an ideal of $R$, then
$$N = \ann_{G}\left(\fa\fc R[x,f]\right)
= \ann_{G}\left((\fa\cap\fc) R[x,f]\right) \in \mathcal{A}(G).$$
Furthermore, if $\fa$ is the ideal in $\mathcal{I}(G)$
corresponding to $L$ and $\fc$ is the ideal in $\mathcal{I}(G/L)$
corresponding to $N/L$, then $\fa\cap\fc$ is the ideal in
$\mathcal{I}(G)$ corresponding to $N$.

\item There is an order-preserving bijection from $\{ N \in
\mathcal{A}(G) : N \supseteq L\}$ to $\mathcal{A}(G/L)$ given by
$N \mapsto N/L$.
\end{enumerate}
\end{prop}

\begin{proof} By Lemma \ref{ga.1}, the
left $R[x,f]$-module $G/L$ is $x$-torsion-free.

(i) Let $g \in N$. Let $i,j \in \nn$ and $r \in (\fb:\fa)$, $u \in
\fa$. Then $ux^i(rx^jg) = ur^{p^i}x^{i+j}g = 0$ because $ ur^{p^i}
\in \fb$ and $\fb x^{i+j}$ annihilates $N$. This is true for all
$i \in \nn$ and $u \in \fa$. Therefore $ rx^jg \in \ann_G\left(\fa
R[x,f]\right) = L. $ Since this is true for all $j \in \nn$ and $r
\in (\fb:\fa)$, we see that $N/L \subseteq
\ann_{G/L}\left((\fb:\fa) R[x,f]\right)$.

Now suppose that $g \in G$ is such that $g + L \in
\ann_{G/L}\left((\fb:\fa) R[x,f]\right)$. Let $r \in \fb$ and $i
\in \nn$. Then $r \in (\fb:\fa)$ and so $rx^{i+1}$ annihilates $g
+ L \in G/L$.  Hence $rx^{i+1}g \in L$. Since $\fb \subseteq \fa$,
we see that $r^{p-1}rx^{i+1}g = 0$, so that $xrx^ig = 0$. As $G$
is $x$-torsion-free, it follows that $rx^ig = 0$. As this is true
for all $r \in \fb$ and $i \in \nn$, we see that $ g \in
\ann_G\left(\fb R[x,f]\right) = N. $ Hence  $N/L =
\ann_{G/L}\left((\fb:\fa) R[x,f]\right)$.

To prove the final claim, we have to show that $\grann_{R[x,f]}(N/L)
= \bigoplus _{n\in\nn} (\fb:\fa) x^n$, given that $\grann_{R[x,f]}N
= \bigoplus _{n\in\nn} \fb x^n$. In view of the preceding
paragraph, it remains only to show that $$\grann_{R[x,f]}(N/L) \subseteq
(\fb:\fa)R[x,f].$$

Let $r \in R$ be such that $rx^i \in \grann_{R[x,f]}(N/L)$ for all
$i \in \nn$. Let $g \in N$. Then $rx^ig \in L$ for all $i \in
\nn$, and so $\fa rx^ig = 0$ for all $i \in \nn$. As this is true
for all $g \in N$ and for all $i \in \nn$, it follows that $\fa r
\subseteq \left(\grann_{R[x,f]}N\right) \cap R = \fb$. Hence $r
\in (\fb:\fa)$.

(ii) Let $g \in N$. Then $ux^ig \in L$ for all $u \in \fc$ and $i
\in \nn$, and so $rux^ig = 0$ for all $r \in \fa$, $u \in \fc$ and
$i \in \nn$. Hence $N \subseteq
\ann_G\left(\bigoplus_{n\in\nn}\fa\fc x^n\right) =
\ann_G\left(\fa\fc R[x,f]\right)$.

Now let $g \in \ann_G\left(\fa\fc R[x,f]\right)$. Then, for all $r
\in \fa$, $u \in \fc$ and $i,j \in \nn$, we have $rx^i(ux^jg) =
ru^{p^i}x^{i+j}g = 0$, and so $ux^jg \in L$ for all $u \in \fc$
and $j \in \nn$. Hence $g + L \in
\ann_{G/L}\left(\bigoplus_{n\in\nn}\fc x^n\right) = N/L$, and $g
\in N$. It follows that $ N = \ann_G\left(\fa\fc R[x,f]\right). $
Also, by Lemma \ref{ga.1a}, we have
$$\ann_G\left(\fa\fc R[x,f]\right) = \ann_G\left((\fa\cap\fc) R[x,f]\right),$$
because $\fa \fc$ and $\fa \cap \fc$ have the same radical.

To prove the final claim, we have to show that $\grann_{R[x,f]}N =
(\fa\cap\fc)R[x,f]$, given that
$$\grann_{R[x,f]}(N/L) = \fc R[x,f] \quad \mbox{~and~} \quad
\grann_{R[x,f]}(L) = \fa R[x,f].$$ In view of the
preceding paragraph, it remains only to show that $\grann_{R[x,f]}N
\subseteq (\fa\cap\fc)R[x,f].$ However, this is clear, because
$
\grann_{R[x,f]}N \subseteq \grann_{R[x,f]}L \cap \grann_{R[x,f]}(N/L).
$

(iii) This is now immediate from parts (i) and (ii).
\end{proof}

\begin{rmk}
\label{ga.3r} It follows from Proposition \ref{ga.3}(ii) (and with
the hypotheses and notation thereof) that, if $\fa$ is an ideal of
$R$ and $L := \ann_G\left(\fa R[x,f]\right)$, then $
\ann_{G/L}\left(\fa R[x,f]\right) = 0$.
\end{rmk}

Because the special $R$-ideals introduced in Definition \ref{nt.6} are
radical, the following lemma will be very useful.

\begin{lem}
\label{ga.2} Let $\fa$ and $\fb$ be proper radical ideals of $R$, and let
their (unique) minimal primary decompositions be
$$
\fa = \fr_1 \cap \ldots \cap \fr_k \cap \fp_1 \cap \ldots \cap \fp_t \cap
\fp_1' \cap \ldots \cap \fp_u'
$$
and
$$
\fb = \fr_1 \cap \ldots \cap \fr_k \cap \fq_1 \cap \ldots \cap \fq_v \cap
\fq_1' \cap \ldots \cap \fq_w',
$$
where the notation is such that
$$
\left\{\fp_1, \ldots, \fp_t,
\fp_1', \ldots, \fp_u'\right\} \cap \left\{\fq_1, \ldots, \fq_v,
\fq_1', \ldots, \fq_w'\right\} = \emptyset,
$$
and such that none of $\fp_1, \ldots, \fp_t$ contains an
associated prime of $\fb$, each of $\fp_1', \ldots, \fp_u'$
contains an associated prime of $\fb$, none of $\fq_1, \ldots,
\fq_v$ contains an associated prime of $\fa$, and each of $\fq_1',
\ldots, \fq_w'$ contains an associated prime of $\fa$. (Note that
some, but not all, of the integers $k$, $t$ and $u$ might be zero;
a similar comment applies to the primary decomposition of $\fb$.)
Then
\begin{enumerate}
\item $\fa \cap \fb = \fr_1 \cap \ldots \cap \fr_k \cap \fp_1 \cap
\ldots \cap \fp_t \cap \fq_1 \cap \ldots \cap \fq_v$ is the
minimal primary decomposition; \item if $\fa \not\subseteq \fb$,
the equation
$ (\fb : \fa) = \fq_1 \cap \ldots \cap \fq_v$
gives the minimal primary decomposition.
\end{enumerate}
\end{lem}

\begin{proof} (i) Each of $\fp_1', \ldots, \fp_u'$ must contain one of
$\fq_1, \ldots, \fq_v$; likewise, each of $\fq_1', \ldots, \fq_w'$
must contain one of $\fp_1, \ldots, \fp_t$ . The claim then
follows easily.

(ii) Since $ (\fb : \fa) = (\fb \cap \fa : \fa)$, it is clear from part (i) that
$$
\fq_1 \cap \ldots \cap \fq_v \subseteq (\fb : \fa).
$$
Now let $r \in (\fb : \fa) = (\fb \cap \fa : \fa)$.
Then, for each $i = 1, \ldots, v$,
we have $r\fa \subseteq \fq_i$, whereas $\fa \not\subseteq \fq_i$;
hence $r \in \fq_i$
because $\fq_i$ is prime.
\end{proof}

\begin{thm}
\label{ga.4} Suppose that $G$ is $x$-torsion-free. Let $N :=
\ann_G\left(\fb R[x,f]\right) \in \mathcal{A}(G)$, where the ideal
$\fb \in \mathcal{I}(G)$ corresponds to $N$. Assume that $N \neq
0$, and let $\fb = \fp_1 \cap \ldots \cap \fp_t$ be the minimal
primary decomposition of the (radical) ideal $\fb$.

Suppose that $t > 1$, and consider any partition $\{1,\ldots,t\} =
U \cup V$, where $U$ and $V$ are two non-empty disjoint sets. Set
$ \fa = \bigcap_{i\in U} \fp_i$ and $\fc = \bigcap_{i\in V}
\fp_i$. Let $L := \ann_G\left(\fa R[x,f]\right) \in
\mathcal{A}(G)$. Then

\begin{enumerate}
\item $0 \subset L \subset N$ (the symbol `$\subset$' is reserved
to denote strict inclusion); \item $N/L = \ann_{G/L}\left(\fc
R[x,f]\right) \in \mathcal{A}(G/L)$ with corresponding ideal $\fc
\in \mathcal{I}(G/L)$; and \item $\grann_{R[x,f]}L = \fa R[x,f]$,
so that $\fa \in \mathcal{I}(G)$ corresponds to $L$.
\end{enumerate}
\end{thm}

\begin{proof} (i) It is clear that $ L \subseteq N$.
Suppose that $L = 0$ and seek a contradiction. Let $g \in N$. Let
$i,j \in \nn$ and $r \in \fc$, $u \in \fa$. Then $ux^i(rx^jg) =
ur^{p^i}x^{i+j}g = 0$ because $ ur^{p^i} \in \fb$ and $\fb
x^{i+j}$ annihilates $N$. This is true for all $i \in \nn$ and $u
\in \fa$. Therefore $ rx^jg \in \ann_G\left(\bigoplus_{n\in\nn}\fa
x^n\right) = L = 0. $ It follows that $\bigoplus_{n\in\nn}\fc x^n
\subseteq \grann_{R[x,f]}N = \bigoplus _{n\in\nn} \fb x^n$, so that
$\fc \subseteq \fb$. But $\fb \subseteq \fc$, and so $\fc = \fb$.
However, this contradicts the fact that $\fb = \fp_1 \cap \ldots
\cap \fp_t$ is the unique minimal primary decomposition of $\fb$.
Therefore $L \neq 0$.

Now suppose that $L = N$ and again seek a contradiction. Then
$\bigoplus_{n\in\nn}\fa x^n \subseteq
\grann_{R[x,f]}N = \bigoplus _{n\in\nn} \fb x^n$, so that $\fa \subseteq \fb$. But
$\fb \subseteq \fa$, and so $\fa = \fb$, and this again leads to a contradiction.
Therefore $L \neq N$.

(ii) Since $\fb \subseteq \fa$, it is immediate from Proposition
\ref{ga.3}(i) that $N/L = \ann_{G/L}\left((\fb:\fa)R[x,f]\right)
\in \mathcal{A}(G/L)$ and that the ideal $(\fb:\fa) \in
\mathcal{I}(G/L)$ corresponds to $N/L$. However, it follows from
Lemma \ref{ga.2}(ii) that $(\fb:\fa) = \bigcap_{i\in V} \fp_i =
\fc$.

(iii) Let $\fd \in \mathcal{I}(G)$ correspond to $L$. Note that
$\fa = \bigcap_{i\in U} \fp_i \subseteq \fd$. By Proposition
\ref{ga.3}(i), the ideal in $\mathcal{I}(G/L)$ corresponding to
$N/L$ is $(\fb:\fd)$. Therefore, by part (ii), we have $(\fb:\fd)
= \fc$. But, by Proposition \ref{ga.3}(ii), the ideal in
$\mathcal{I}(G)$ corresponding to $N$ is $\fd\cap\fc$. Therefore
$\fb = \fd \cap \fc$, and so $\fd \neq R$.

Now $\fd$ is a radical ideal of $R$. By Lemma \ref{ga.2}(i), each
$\fp_j$, for $j \in U$, is an associated prime of $\fd$. Hence
$\fd \subseteq \bigcap_{j\in U} \fp_j = \fa$. But we already know
that $\fa \subseteq \fd$, and so $\fd = \fa$.
\end{proof}

\begin{cor}
\label{ga.10} Suppose
that $G$ is $x$-torsion-free. Then the set of $G$-special
$R$-ideals is precisely the set of all finite intersections of
prime $G$-special $R$-ideals (provided one includes the
empty intersection, $R$, which corresponds to the zero special annihilator
submodule of $G$). In symbols,
$$
\mathcal{I}(G) = \left\{ \fp_1 \cap \ldots \cap \fp_t : t \in \nn
\mbox{~and~} \fp_1, \ldots, \fp_t \in
\mathcal{I}(G)\cap\Spec(R)\right\}.
$$
\end{cor}

\begin{proof} By Corollary \ref{nt.9}, the set $\mathcal{I}(G)$ is
closed under taking intersections. A proper ideal $\fa \in
\mathcal{I}(G)$ is radical and it follows from Theorem \ref{ga.4}
that each (necessarily prime) primary component of $\fa$ also
belongs to $\mathcal{I}(G)$. This is enough to complete the proof.
\end{proof}

\begin{lem}
\label{ga.11} Suppose that $G$ is $x$-torsion-free. Let $\fp$ be a
maximal member of $\mathcal{I}(G) \setminus \{R\}$ with respect to
inclusion, and let $L \in \mathcal{A}(G)$ be the corresponding
special annihilator submodule of $G$. Thus $L$ is a minimal member
of the set of non-zero special annihilator submodules of $G$.

Then $\fp$ is prime, and any non-zero $g\in L$ satisfies
$\grann_{R[x,f]}R[x,f]g = \fp R[x,f]$.
\end{lem}

\begin{proof} It follows from Corollary \ref{ga.10} that $\fp$ is prime.

Since $R[x,f]g$ is a non-zero $R[x,f]$-submodule of $L$, there is
a proper radical ideal $\fa \in \mathcal{I}(G)$ such that
$$
\fa R[x,f] = \grann_{R[x,f]}R[x,f]g \supseteq \grann_{R[x,f]}L = \fp
R[f,x].
$$
Since $\fp$ is a maximal member of $\mathcal{I}(G) \setminus
\{R\}$, we must have $\fa = \fp$.
\end{proof}

Our next major aim is to show that, in the situation of Corollary
\ref{ga.10}, the set $\mathcal{I}(G)$ is finite if $G$ has the
property that, for each special annihilator submodule $L$ of $G$
(including $0 = \ann_GR[x,f]$), the $x$-torsion-free residue class
module $G/L$ (see Lemma \ref{ga.1}) does not contain, as an
$R[x,f]$-submodule, an infinite direct sum of non-zero special
annihilator submodules of $G/L$. This may seem rather a
complicated hypothesis, and so we point out now that it is
satisfied if $G$ is a Noetherian or Artinian left $R[x,f]$-module,
and therefore if $G$ is a Noetherian or Artinian $R$-module. These
ideas will be applied, later in the paper, to an example in which
$G$ is Artinian as an $R$-module.

The following lemma will be helpful in an inductive argument in
the proof of Theorem \ref{ga.12}.

\begin{lem}
\label{ga.12p} Suppose that $G$ is $x$-torsion-free, and that the
set $\mathcal{I}(G) \setminus \{R\}$ is non-empty and has finitely
many maximal members: suppose that there are $n$ of these and
denote them by $\fp_1, \ldots, \fp_n$. (The ideals $\fp_1, \ldots,
\fp_n$ are prime, by \/ {\rm \ref{ga.11}}.) Let $L :=
\ann_G\left(\fp_1 \cap \cdots \cap \fp_n\right)R[x,f]$. Then the
left $R[x,f]$-module $G/L$ is $x$-torsion-free, and
$$
\mathcal{I}(G/L) \cap \Spec (R) = \mathcal{I}(G) \cap \Spec (R)
\setminus \{\fp_1, \ldots, \fp_n\}.
$$
\end{lem}

\begin{proof} Note
that $\bigcap_{i=1}^n \fp_i \in \mathcal{I}(G)$, by Corollary
\ref{nt.9}. Therefore $\grann_{R[x,f]}L = \left(\bigcap_{i=1}^n
\fp_i\right)R[x,f]$ and $L$ corresponds to $\bigcap_{i=1}^n
\fp_i$.

That $G/L$ is $x$-torsion-free follows from Lemma \ref{ga.1}. By
Proposition \ref{ga.3}(iii),
$$
\mathcal{A}(G/L) = \left\{ N/L : N \in \mathcal{A}(G) \mbox{~and~}
L \subseteq N \right\}.
$$
Let $N \in \mathcal{A}(G)$ with $L \subset N$, and let $\fb \in
\mathcal{I}(G)$ correspond to $N$. Note that $\fb \subset
\bigcap_{i=1}^n \fp_i$, and that no associated prime of $\fb$ can
contain properly any of $\fp_1, \ldots, \fp_n$. Therefore the
minimal primary decomposition of the radical ideal $\fb$ will have
the form
$$
\fb = \left( {\textstyle \bigcap_{i\in I}} \fp_i \right) \cap
\fq_1 \cap \ldots \cap \fq_v,
$$
where $I$ is some (possibly empty) subset of $\{1, \ldots, n\}$
and none of $\fq_1, \ldots, \fq_v$ contains any of $\fp_1, \ldots,
\fp_n$. Note that $\fq_1, \ldots, \fq_v$ must all belong to
$\mathcal{I}(G) \cap \Spec (R) \setminus \{\fp_1, \ldots,
\fp_n\}$. Proposition \ref{ga.3}(i), this time used in conjunction
with Lemma \ref{ga.2}(ii), now shows that $N/L \in
\mathcal{A}(G/L)$ and the ideal of $\mathcal{I}(G/L)$
corresponding to $N/L$ is
$$
\left(\fb : \fp_1 \cap \cdots \cap \fp_n\right) = \fq_1 \cap
\ldots \cap \fq_v.
$$
Note also that, if $\fq \in \mathcal{I}(G) \cap \Spec (R)
\setminus \{\fp_1, \ldots, \fp_n\}$ and
$$
J := \left\{ j \in \{1, \ldots, n\} : \fp_j \not\supset \fq
\right\},
$$
then $\fc :=  \left(\bigcap_{j\in J} \fp_j \right) \cap \fq \in
\mathcal{I}(G)$ and $\fc \subset \bigcap_{i=1}^n\fp_i$. It now
follows from Corollary \ref{ga.10} that
$$
\mathcal{I}(G/L) \cap \Spec (R) = \mathcal{I}(G) \cap \Spec (R)
\setminus \{\fp_1, \ldots, \fp_n\},
$$
as required.
\end{proof}

\begin{thm}
\label{ga.12} Suppose
that $G$ is $x$-torsion-free. Assume that $G$ has the property
that, for each special annihilator submodule $L$ of $G$ (including $0 =
\ann_GR[x,f]$), the $x$-torsion-free residue class module $G/L$
does not contain, as an $R[x,f]$-submodule, an infinite direct sum
of non-zero special annihilator submodules of $G/L$.

Then the set $\mathcal{I}(G)$ of $G$-special $R$-ideals is
finite.
\end{thm}

\begin{proof} By Corollary \ref{ga.10}, it is enough for us to show that
the set $\mathcal{I}(G)\cap \Spec (R)$ is finite; we may suppose
that the latter set is not empty, so that it has maximal members
with respect to inclusion. In the first part of the proof, we show
that $\mathcal{I}(G)\cap \Spec (R)$ has only finitely many such
maximal members.

Let $\left(\fp_{\lambda}\right)_{\lambda \in \Lambda}$ be a
labelling of the set of maximal members of $\mathcal{I}(G)\cap
\Spec (R)$, arranged so that $\fp_{\lambda} \neq \fp_{\mu}$
whenever $\lambda$ and $\mu$ are different elements of $\Lambda$.
For each $\lambda \in \Lambda$, let $S_{\lambda}$ be the member of
$\mathcal{A}(G)$ corresponding to $\fp_{\lambda}$.

Consider $\lambda, \mu \in \Lambda$ with $\lambda \neq \mu$. By
Lemma \ref{ga.11}, a non-zero $g \in S_{\lambda} \cap S_{\mu}$
would have to satisfy $\grann_{R[x,f]}R[x,f]g = \fp_{\lambda} R[x,f]
= \fp_{\mu} R[x,f]$. Since $\fp_{\lambda} \neq \fp_{\mu}$, this is
impossible. Therefore $S_{\lambda} \cap S_{\mu} = 0$ and the sum
$S_{\lambda} + S_{\mu}$ is direct.

Suppose, inductively, that $n \in \N$ and we have shown that,
whenever $\lambda_1, \ldots, \lambda_n$ are $n$ distinct members
of $\Lambda$, then the sum $\sum_{i=1}^nS_{\lambda_i}$ is direct.
We can now use Lemma \ref{ga.11} to see that, if $g_i \in
S_{\lambda_i}$ for $i = 1, \ldots, n$, then
$$
\grann_{R[x,f]}R[x,f](g_1 + \cdots + g_n) =
\bigcap_{\stackrel{\scriptstyle i=1}{g_i \neq
0}}^n\fp_{\lambda_i}R[x,f],
$$
and then to deduce that, for $\lambda_{n+1}
\in \Lambda \setminus\{\lambda_1, \ldots,
\lambda_n\}$, we must have
$
\left( \bigoplus_{i=1}^n S_{\lambda_i} \right) \bigcap S_{\lambda_{n+1}} = 0,
$
so that the sum $S_{\lambda_1} + \cdots + S_{\lambda_n} + S_{\lambda_{n+1}}$
is direct.

It follows that the sum $\sum_{\lambda \in \Lambda}S_{\lambda}$ is
direct; since each $S_{\lambda}$ is non-zero, the hypothesis about
$G/0$ (that is, about $G$) ensures that $\Lambda$ is finite.

We have thus shown that $\mathcal{I}(G)\cap \Spec (R)$ has only
finitely many maximal members. Note that $\max\{ \height \fp : \fp
\mbox{~is a maximal member of~} \mathcal{I}(G)\cap \Spec (R) \}$
is an upper bound for the lengths of chains
$$
\fp_0 \subset \fp_1 \subset \cdots \subset \fp_w
$$
of prime ideals in $\mathcal{I}(G)\cap \Spec (R)$. We argue by
induction on the maximum $t$ of these lengths. When $t = 0$, all
members of $\mathcal{I}(G)\cap \Spec (R)$ are maximal members of
that set, and so, by the first part of this proof,
$\mathcal{I}(G)\cap \Spec (R)$ is finite. Now suppose that $t >
0$, and that it has been proved that $\mathcal{I}(G)\cap \Spec
(R)$ is finite for smaller values of $t$.

We know that there are only finitely many maximal members of
$\mathcal{I}(G)\cap \Spec (R)$; suppose that there are $n$ of
these and denote them by $\fp_1, \ldots, \fp_n$. Let $L :=
\ann_G\left(\fp_1 \cap \cdots \cap \fp_n\right)R[x,f]$. We can now
use Lemma \ref{ga.12p} to deduce that the left $R[x,f]$-module
$G/L$ is $x$-torsion-free and
$$
\mathcal{I}(G/L) \cap \Spec (R) = \mathcal{I}(G) \cap \Spec (R)
\setminus \{\fp_1, \ldots, \fp_n\}.
$$
It follows from this and Proposition \ref{ga.3}(ii) that the
inductive hypothesis can be applied to $G/L$, and so we can deduce
that the set $$\mathcal{I}(G) \cap \Spec (R) \setminus \{\fp_1,
\ldots, \fp_n\}$$ is finite. Hence $\mathcal{I}(G) \cap \Spec (R)$
is a finite set and the inductive step is complete.
\end{proof}

\begin{cor}
\label{ga.12c} Suppose that the left $R[x,f]$-module $G$ is
$x$-torsion-free and either Artinian or Noetherian as an
$R$-module. Then the set $\mathcal{I}(G)$ of $G$-special
$R$-ideals is finite.
\end{cor}

\begin{thm}
\label{ga.13} Suppose that $G$ is $x$-torsion-free and that the
set $\mathcal{I}(G)$ of $G$-special $R$-ideals is finite.
Then there exists a (uniquely determined) ideal $\fb \in
\mathcal{I}(G)$ with the properties that $\height \fb \geq 1$ (the
improper ideal $R$ is considered to have infinite height) and $\fb
\subset \fc$ for every other ideal $\fc \in \mathcal{I}(G)$ with
$\height \fc \geq 1$. Furthermore, for $g \in G$, the following
statements are equivalent:

\begin{enumerate}
\item $g$ is annihilated by $\fb R[x,f] = \bigoplus_{n\in\nn}\fb x^n$;
\item there exists $c \in R^{\circ}\cap \fb$ such that
$cx^ng = 0$ for all $n \gg 0$;
\item there exists $c \in R^{\circ}$ such that $cx^ng = 0$ for all $n \gg 0$.
\end{enumerate}
\end{thm}

\begin{proof} By Corollary \ref{ga.10}, we have
$$
\mathcal{I}(G) = \left\{ \fp_1 \cap \ldots \cap \fp_t : t \in \nn
\mbox{~and~} \fp_1, \ldots, \fp_t \in
\mathcal{I}(G)\cap\Spec(R)\right\}.
$$
Since $\mathcal{I}(G)$ is finite, it is immediate that
$$
\fb := \bigcap_{\stackrel{\scriptstyle \fp \in
\mathcal{I}(G)\cap\Spec(R)}{\height \fp \geq 1}}\fp
$$
is the smallest ideal in $\mathcal{I}(G)$ of height greater than
$0$. Since $\height \fb \geq 1$, so that there exists $c \in \fb
\cap R^{\circ}$ by prime avoidance, it is clear that (i)
$\Rightarrow$ (ii) and (ii) $\Rightarrow$ (iii).

(iii) $\Rightarrow$ (i) Let $n_0 \in \nn$ and $c \in R^{\circ}$ be
such that $cx^ng = 0$ for all $n \geq n_0$. Then, for all $j \in \nn$,
we have $x^{n_0}cx^jg = c^{p^{n_0}}x^{n_0 + j}g = 0$, so that $cx^jg = 0$
because $G$ is $x$-torsion-free.

Therefore $g \in \ann_G(RcR[x,f])$. Now $\ann_G(RcR[x,f])\in
\mathcal{A}(G)$: let $\fa \in \mathcal{I}(G)$ be the corresponding
$G$-special $R$-ideal. Since $c \in \fa$, we must have
$\height \fa \geq 1$. Therefore $\fb \subseteq \fa$, by definition
of $\fb$, and so
$$
g \in \ann_G(RcR[x,f]) = \ann_G(\fa R[x,f]) \subseteq \ann_G(\fb
R[x,f]).
$$
\end{proof}

Corollary \ref{ga.12c} and Theorem \ref{ga.13} give hints about
how this work will be exploited, in Section \ref{tc} below, to
obtain results in the theory of tight closure. The aim is to apply
Corollary \ref{ga.12c} and Theorem \ref{ga.13} to
$H^d_{\fm}(R)/\Gamma_x(H^d_{\fm}(R))$, where $(R,\fm)$ is a local
ring of dimension $d > 0$; the local cohomology module
$H^d_{\fm}(R)$, which is well known to be Artinian as an
$R$-module, carries a natural structure as a left $R[x,f]$-module.
The passage between $H^d_{\fm}(R)$ and its $x$-torsion-free
residue class $R[x,f]$-module
$H^d_{\fm}(R)/\Gamma_x(H^d_{\fm}(R))$ is facilitated by the
following extension, due to G. Lyubeznik, of a result of R.
Hartshorne and R. Speiser. It shows that, when $R$ is local, an
$x$-torsion left $R[x,f]$-module which is Artinian (that is,
`cofinite' in the terminology of Hartshorne and Speiser) as an
$R$-module exhibits a certain uniformity of behaviour.

\begin{thm} [G. Lyubeznik {\cite[Proposition 4.4]{Lyube97}}]
\label{hs.4}  {\rm (Compare Hartshorne--Speiser \cite[Proposition
1.11]{HarSpe77}.)} Suppose that $(R,\fm)$ is local, and let $H$ be
a left $R[x,f]$-module which is Artinian as an $R$-module. Then
there exists $e \in \N_0$ such that $x^e\Gamma_x(H) = 0$.
\end{thm}

Hartshorne and Speiser first proved this result in the
case where $R$ is local and contains its residue field which is perfect.
Lyubeznik applied his theory of $F$-modules to obtain the result
without restriction on the local ring $R$ of characteristic $p$.

\begin{defi}
\label{hslno} Suppose that $(R,\fm)$ is local, and let $H$ be a
left $R[x,f]$-module which is Artinian as an $R$-module.  By the
Hartshorne--Speiser--Lyubeznik Theorem \ref{hs.4}, there exists $e
\in \nn$ such that $x^e\Gamma_x(H) = 0$: we call the smallest such
$e$ the {\em Hartshorne--Speiser--Lyubeznik number\/}, or {\em
HSL-number\/} for short, of $H$.
\end{defi}

It will be helpful to have available an extension of this idea.

\begin{defi}
\label{ga.14} We say that the left $R[x,f]$-module $H$ {\em admits
an HSL-number\/} if there exists $e \in \nn$ such that $x^e\Gamma_x(H) = 0$;
then we call the smallest such $e$ the {\em HSL-number\/} of $H$.
\end{defi}

We have seen above in \ref{hs.4} and \ref{hslno} that if $H$ is Artinian
as an $R$-module, then it admits an HSL-number.  Note also that
if $H$ is Noetherian
as an $R$-module, then it admits an HSL-number, because $\Gamma_x(H)$ is an
$R[x,f]$-submodule of $H$, and so is an $R$-submodule and therefore
finitely generated.

\begin{cor}
\label{ga.15} Suppose that the left $R[x,f]$-module $H$ admits an
HSL-number $m_0$, and that the $x$-torsion-free left
$R[x,f]$-module $G := H/\Gamma_x(H)$ has only finitely many
$G$-special $R$-ideals. Let $\fb$ be the smallest ideal in
$\mathcal{I}(G)$ of positive height (see\/ {\rm \ref{ga.13}}). For
$h \in H$, the following statements are equivalent:

\begin{enumerate}
\item $h$ is annihilated by $\bigoplus_{n\geq m_0}\fb^{[p^{m_0}]} x^n$;
\item there exists $c \in R^{\circ}\cap \fb$ such that
$cx^nh = 0$ for all $n \geq m_0$;
\item there exists $c \in R^{\circ}\cap \fb$ such that
$cx^nh = 0$ for all $n \gg 0$;
\item there exists $c \in R^{\circ}$ such that $cx^nh = 0$ for all $n \gg 0$.
\end{enumerate}
\end{cor}

\begin{proof} Since $\fb
\cap R^{\circ} \neq 0$ by prime avoidance, it is clear that (i)
$\Rightarrow$ (ii), (ii)
$\Rightarrow$ (iii) and (iii) $\Rightarrow$ (iv).

(iv) $\Rightarrow$ (i) Since $cx^n(h + \Gamma_x(H)) = 0$ in $G$
for all $n \gg 0$, it follows from Theorem \ref{ga.13} that $h + \Gamma_x(H)$
is annihilated by $\fb R[x,f]$. Therefore, for all $r \in \fb$ and $j \in \nn$,
we have $rx^j(h + \Gamma_x(H)) = 0$, so that $rx^jh \in \Gamma_x(H)$ and
$r^{p^{m_0}}x^{m_0+j}h = x^{m_0}rx^jh = 0$. Therefore $h \in \ann_H
\left(\bigoplus_{n\geq m_0}\fb^{[p^{m_0}]} x^n\right)$.
\end{proof}

\section{\sc Applications to tight closure}
\label{tc}

The aim of this section is to apply results from Section \ref{ga}
to the theory of tight closure in the local ring $(R,\fm)$ of
dimension $d > 0$. As was mentioned in Section \ref{ga}, we shall
be concerned with the top local cohomology module $H^d_{\fm}(R)$,
which has a natural structure as a left $R[x,f]$-module, and its
$x$-torsion-free residue class module
$H^d_{\fm}(R)/\Gamma_x(H^d_{\fm}(R))$. The (well-known) left
$R[x,f]$-module structure carried by $H^d_{\fm}(R)$ is described
in detail in \cite[2.1 and 2.3]{KS}.

\begin{rmd}
\label{tc.1} Suppose that $(R,\fm)$ is a local ring of dimension $d > 0$.
The above-mentioned natural left $R[x,f]$-module structure carried
by $H^d_{\fm}(R)$ is independent of any choice of a system of parameters for
$R$.  However, if one does choose a system of parameters $a_1, \ldots, a_d$ for
$R$, then one can obtain a quite concrete
representation of the local cohomology module $H^d_{\fm}(R)$ and, through this,
an explicit formula for the effect of multiplication by the indeterminate
$x \in R[x,f]$ on an element of $H^d_{\fm}(R)$.

Denote by $a_1, \ldots, a_d$ a system of parameters for $R$.

\begin{enumerate}
\item Represent
$H^d_{\fm}(R)$ as the $d$-th cohomology module of the \u{C}ech
complex of $R$ with respect to $a_1, \ldots, a_d$, that is, as the
residue class module of $R_{a_1 \ldots a_d}$ modulo the image,
under the \u{C}ech complex `differentiation' map, of
$\bigoplus_{i=1}^dR_{a_1 \ldots a_{i-1}a_{i+1}\ldots a_d}$. See
\cite[\S 5.1]{LC}. We use `$\left[\phantom{=} \right]$' to denote natural
images of elements of $R_{a_1\ldots a_d}$ in this residue class
module. Note that, for $i \in \{1, \ldots, d\}$, we have
$$
\left[\frac{a_i^k}{(a_1 \ldots a_d)^k}\right] = 0 \quad \mbox{~for
all~} k \in \nn.
$$

Denote the product $a_1 \ldots a_d$ by $a$. A typical element of
$H^d_{\fm}(R)$ can be represented as $ \left[r/a^j\right]$ for
some $r \in R$ and $j \in \nn$; moreover, for $r, r_1 \in R$ and
$j, j_1 \in \nn$, we have $ \left[r/a^j\right] =
\left[r_1/a^{j_1}\right] $ if and only if there exists $k \in \nn$
such that $k \geq \max\{j,j_1\}$ and $ a^{k-j}r - a^{k-j_1}r_1 \in
(a_1^k, \ldots, a_d^k)R. $ In particular, if $a_1, \ldots, a_d$
form an $R$-sequence (that is, if $R$ is Cohen--Macaulay), then $
\left[r/a^j\right] = 0$ if and only if $r \in (a_1^j, \ldots,
a_d^j)R$, by \cite[Theorem 3.2]{O'Car83}, for example.

\item The left $R[x,f]$-module structure on $H^d_{\fm}(R)$ is such
that
$$
x\left[\frac{r}{(a_1\ldots a_d)^j}\right] =
\left[\frac{r^p}{(a_1\ldots a_d)^{jp}}\right] \quad \mbox{~for all~} r \in R
\mbox{~and~} j \in \nn.
$$
The reader might like to consult \cite[2.3]{KS} for more details, and should
in any case note that this left $R[x,f]$-module structure does not depend on the
choice of system of parameters $a_1, \ldots, a_d$.
\end{enumerate}
\end{rmd}

\begin{rmk}
\label{tc.1r} Let the situation and notation be as in \ref{tc.1}.
Here we relate the left $R[x,f]$-module structure on $H :=
H^d_{\fm}(R)$ described in \ref{tc.1} to the tight closure in $H$
of its zero submodule. See \cite[Definition (8.2)]{HocHun90} for
the definition of the tight closure in an $R$-module of one of its
submodules. Let $n \in \nn$.

\begin{enumerate}
\item The $n$-th component $Rx^n$ of $R[x,f]$ is isomorphic, as an
$(R,R)$-bimodule, to $R$ considered as a left $R$-module in the
natural way and as a right $R$-module via $f^n$, the $n$-th power
of the Frobenius ring homomorphism. Let $L$ be a submodule of the
$R$-module $M$. It follows that an element $m \in M$ belongs to
$L^*_M$, the {\em tight closure of $L$ in $M$\/}, if and only if
there exists $c \in R^\circ$ such that $cx^n \otimes m$ belongs,
for all $n \gg 0$, to the image of $R[x,f]\otimes_R L$ in
$R[x,f]\otimes_R M$ under the map induced by inclusion.

\item Let $S$ be a multiplicatively closed subset of $R$. It is straightforward
to check that there is an isomorphism of $R$-modules
$$\gamma_n: Rx^n \otimes_R S^{-1}R \stackrel{\cong}{\lra} S^{-1}R$$
for which $\gamma_n ( bx^n \otimes (r/s) ) = br^{p^n}/s^{p^n}$ for
all $b,r \in R$ and $s \in S$; the inverse of $\gamma_n$ satisfies
$(\gamma_n)^{-1} (r/s) = rs^{p^n-1}x^n \otimes (1/s)$ for all $r \in
R$ and $s \in S$.

\item Now represent $H := H^d_{\fm}(R)$ as the $d$-th cohomology
module of the \u{C}ech complex of $R$ with respect to the system
of parameters $a_1, \ldots, a_d$, as in \ref{tc.1}(i). We can use
isomorphisms like that described in part (ii), together with the
right exactness of tensor product, to see that (when we think of
$H$ simply as an $R$-module) there is an isomorphism of
$R$-modules $\delta_n: Rx^n \otimes_R H \stackrel{\cong}{\lra} H$
for which
$$
\delta_n \left( bx^n \otimes \left[\frac{r}{(a_1\ldots
a_d)^j}\right] \right) = \left[\frac{br^{p^n}}{(a_1\ldots
a_d)^{jp^n}}\right] \quad \mbox{~for all~} b,r \in R \mbox{~and~}
j \in \nn.
$$
Thus, in terms of the natural left $R[x,f]$-module structure on
$H$, we have $\delta_n \left( bx^n \otimes h\right) = bx^nh$ for
all $b \in R$ and $h \in H$.

\item It thus follows that, for $h \in H$, we have $h
\in 0^*_{H}$ if and only if there exists $c \in R^{\circ}$ such
that $cx^nh = 0$ for all $n \gg 0$.

\item Observe that $\Gamma_x(H)
\subseteq 0^*_H$.

\item Suppose that $(R,\fm)$ is Cohen--Macaulay, and use the
notation of part (iii) again; write $a := a_1\ldots a_d$. Let $r
\in R$, $j \in \nn$, and let $h := \left[r/a^j\right]$ in $H$. It
follows from \ref{tc.1}(i) that $r \in ((a_1^j, \ldots,
a_d^j)R)^*$ if and only if there exists $c \in R^{\circ}$ such
that $cx^nh = 0$ for all $n \gg 0$. Thus, by part (iv) above, $r
\in ((a_1^j, \ldots, a_d^j)R)^*$ if and only if $h \in 0^*_{H}$.

\item Let the situation and notation be as in part (vi) above.
Then the $R$-homomorphism $\nu_j : R/(a_1^j, \ldots, a_d^j)R \lra
H$ for which $\nu_j(r' + (a_1^j, \ldots, a_d^j)R) =
\left[r'/a^j\right]$ for all $r'\in R$ is a monomorphism (by
\ref{tc.1}(i)). Furthermore, the induced homogeneous
$R[x,f]$-homomorphism
$$
R[x,f]\otimes_R\nu_j : R[x,f]\otimes_R\left(R/(a_1^j, \ldots,
a_d^j)R\right) \lra R[x,f]\otimes_RH
$$
of graded left $R[x,f]$-modules is also a monomorphism: this is
because a homogeneous element of $\Ker
\left(R[x,f]\otimes_R\nu_j\right)$ must have the form
$r'x^k\otimes(1 + (a_1^j, \ldots, a_d^j)R)$ for some $r' \in R$
and $k \in \nn$; since $r'x^k\otimes\left[1/a^j\right] = 0$, it
follows from \ref{tc.1r}(iii) and \ref{tc.1}(i) that $r' \in
(a_1^{jp^k}, \ldots, a_d^{jp^k})R$, so that $r'x^k\otimes(1 +
(a_1^j, \ldots, a_d^j)R) = 0$.
\end{enumerate}
\end{rmk}

\begin{lem}
\label{tc.1s} Suppose that $(R,\fm)$ is a local ring of dimension $d > 0$;
set $H := H^d_{\fm}(R)$ and $G := H/\Gamma_x(H)$. Let $h \in H$. Then the
following statements are equivalent:

\begin{enumerate}
\item $h \in 0^*_{H}$;
\item $h + \Gamma_x(H) \in 0^*_{G}$;
\item there exists $c \in R^{\circ}$ such
that $cx^n(h+\Gamma_x(H)) = 0$ in $G$ for all $n \gg 0$;
\item there exists $c \in R^{\circ}$ such
that $cx^n(h+\Gamma_x(H)) = 0$ in $G$ for all $n \geq 0$.
\end{enumerate}
\end{lem}

\begin{proof} Let
$m_0$ denote the HSL-number of $H$ (see\/ {\rm \ref{hslno}}).

(i) $\Rightarrow$ (ii) This is immediate from the fact that
$0^*_{H} \subseteq (\Gamma_x(H))^*_H$ once it is recalled from
\cite[Remark (8.4)]{HocHun90} that $h + \Gamma_x(H) \in 0^*_G$ if
and only if $h \in (\Gamma_x(H))^*_H$.

(ii) $\Rightarrow$ (iii) Suppose that $h + \Gamma_x(H) \in
0^*_{G}$, so that $h \in (\Gamma_x(H))^*_H$. Under the isomorphism
$\delta_n: Rx^n \otimes_R H \stackrel{\cong}{\lra}H$ of
\ref{tc.1r}(iii) (where $n \in \nn$), the image of $Rx^n \otimes_R
\Gamma_x(H)$ is mapped into $\Gamma_x(H)$. Therefore there exists
$c \in R^{\circ}$ such that $cx^nh \in\Gamma_x(H)$ for all $n \gg
0$, that is, such that $cx^n(h+\Gamma_x(H)) = 0$ in $G$ for all $n
\gg 0$.

(iii) $\Rightarrow$ (iv) Suppose that there exist
$c \in R^{\circ}$ and $n_0 \in \nn$ such
that $cx^n(h+\Gamma_x(H)) = 0$ in $G$ for all $n \geq n_0$. Then,
for all $j \in \nn$, we have
$$
x^{n_0}cx^j(h+\Gamma_x(H)) = c^{p^{n_0}}x^{n_0 + j}(h+\Gamma_x(H)) = 0 \quad
\mbox{~in~} G.
$$
Since $G$ is $x$-torsion-free, we see that $cx^j(h+\Gamma_x(H)) = 0$ for all
$j \geq 0$.

(iv) $\Rightarrow$ (i) Suppose that there exists $c \in R^{\circ}$ such
that $cx^n(h+\Gamma_x(H)) = 0$ in $G$ for all $n \geq 0$. Then $cx^nh \in
\Gamma_x(H)$ for all $n \geq 0$, so that
$x^{m_0}cx^nh = 0$ for all $n \geq 0$. This implies that
$c^{p^{m_0}}x^{m_0 + n}h = 0$ for all $n \geq 0$, so that $h \in 0^*_{H}$ by
\ref{tc.1r}(iv).
\end{proof}

\begin{defi}
\label{tc.1t} The {\em weak parameter test ideal $\sigma'(R)$
of $R$\/} is defined to be
the ideal generated by $0$ and all weak parameter test elements for $R$.  (By a
`weak parameter test element' for $R$ we mean a $p^{w_0}$-weak
parameter test element for $R$ for some $w_0 \in \nn$.)

It is easy to see that each element of $\sigma'(R) \cap R^{\circ}$
is a weak parameter test element for $R$.
\end{defi}

The next theorem is one of the main results of this paper.

\begin{thm}
\label{tc.2} Let $(R,\fm)$ (as in\/ {\rm \ref{nt.1}})
be a Cohen--Macaulay local ring of
dimension $d > 0$.
Set $H := H^d_{\fm}(R)$, a left $R[x,f]$-module which
is Artinian as an $R$-module; let $m_0$ be its HSL-number (see\/ {\rm \ref{hslno}}),
and let $q_0 := p^{m_0}$.

Set $G := H/\Gamma_x(H)$, an $x$-torsion-free left
$R[x,f]$-module. By\/ {\rm \ref{ga.12c}} and\/ {\rm \ref{ga.13}},
there exists a (uniquely determined) smallest ideal $\fb$ of
height at least $1$ in the set $\mathcal{I}(G)$ of
$G$-special $R$-ideals.

Let $c$ be any element of $\fb \cap R^{\circ}$. Then $c^{q_0}$ is a $q_0$-weak
parameter test element for $R$. In particular, $R$ has a $q_0$-weak
parameter test element. In fact, the weak parameter test ideal $\sigma'(R)$ of
$R$ satisfies $\fb^{[q_0]} \subseteq \sigma'(R) \subseteq \fb$.
\end{thm}

\begin{note} It should be noted that, in Theorem \ref{tc.2}, it is
not assumed that $R$ is excellent. There are examples of Gorenstein
local rings of characteristic $p$
which are not excellent: see \cite[p.\ 260]{HMold}.
\end{note}

\begin{proof}  We have to show that, for an arbitrary parameter ideal $\fa$
of $R$ and $r \in \fa^*$, we have $c^{q_0}r^{p^n} \in \fa^{[p^n]}$ for all
$n \geq m_0$. In the first part of the proof, we establish this in the case where
$\fa$ is an ideal $\fq$ generated by a full system of parameters $a_1, \ldots, a_d$
for $R$.

Let $r \in \fq^*$, so that there exists $\widetilde{c} \in R^{\circ}$ such
that $\widetilde{c}r^{p^n} \in \fq^{[p^n]}$ for all $n \gg 0$. Use $a_1,
\ldots, a_d$ in the notation of \ref{tc.1}(i) for $H^d_{\fm}(R) =
H$, and write $a := a_1\ldots a_d$. We have $ \widetilde{c}x^n\left[r/a\right]
= \left[\widetilde{c}r^{p^n}\!/a^{p^n}\right] = 0 $ in $H$ for all $n \gg 0$.
Set $h := \left[r/a\right] \in H$. Thus $\widetilde{c}x^nh = 0$ for all $n \gg 0$.

It therefore follows from Corollary \ref{ga.15} that $h$ is annihilated by $
\bigoplus_{n\geq m_0}\fb^{[p^{m_0}]} x^n$, so that, in particular,
$c^{p^{m_0}}x^nh = 0$ for all $n \geq m_0$. Hence, in $H$,
$$
\left[\frac{c^{q_0}r^{p^{n}}}{(a_1 \ldots a_d)^{p^{
n}}}\right] = c^{p^{m_0}}x^{n}\left[\frac{r}{a_1 \ldots
a_d}\right] = c^{p^{m_0}}x^n h = 0 \quad \mbox{~for all~} n \geq m_0.
$$
Since $R$ is Cohen--Macaulay, we can now deduce from \ref{tc.1}(i) that
$c^{q_0}r^{p^n} \in \fq^{[p^n]}$ for all $n \geq m_0$, as required (for $\fq$).

Now let $\fa$ be an arbitrary parameter ideal of $R$. A proper
ideal in a Cohen--Macaulay local ring is a parameter ideal if and
only if it can be generated by part of a system of parameters. In
view of the first part of this proof, we can, and do, assume that
$\height \fa < d$. There exist a system of parameters $a_1,
\ldots, a_d$ for $R$ and an integer $i \in \{0, \ldots, d-1\}$
such that $\fa = (a_1,\ldots, a_i)R$. Let $r \in \fa^*$. Then, for
each $v \in \N$, we have $r \in ((a_1,\ldots, a_i,a_{i+1}^v,
\ldots, a_d^v)R)^*$, and, since $a_1,\ldots, a_i,a_{i+1}^v,
\ldots, a_d^v$ is a system of parameters for $R$, it follows from
the first part of this proof that
$$
c^{q_0}r^{p^n} \in ((a_1,\ldots,
a_i,a_{i+1}^v, \ldots, a_d^v)R)^{[p^n]} = \left(a_1^{p^n},\ldots,
a_i^{p^n},a_{i+1}^{vp^n}, \ldots, a_d^{vp^n}\right)\!R
\quad \mbox{~for all~} n \geq m_0.
$$
Therefore, for all $n \geq m_0$,
$$
c^{q_0}r^{p^n} \in \bigcap_{v \in \N} \left(a_1^{p^n},\ldots,
a_i^{p^n},a_{i+1}^{vp^n}, \ldots, a_d^{vp^n}\right)\!R \subseteq
\bigcap_{v \in \N} \left(\fa^{[p^{n}]} + \fm^{vp^{n}} \right) =
\fa^{[p^{n}]}
$$
by Krull's Intersection Theorem. This shows that
$c^{q_0}$ is a $q_0$-weak
parameter test element for $R$, so that, since $\fb$ can be generated by
elements in $\fb \cap R^{\circ}$, it follows that $\fb^{[q_0]} \subseteq \sigma'(R)$.

Now let $c \in \sigma'(R)\cap R^{\circ}$; we suppose that $c \not\in \fb$ and
seek a contradiction. Thus $\fb \subset \fb + Rc$. Let $L : =
\ann_G(\fb R[x,f])$ and $L' := \ann_G((\fb + Rc)R[x,f])$, two
special annihilator submodules of the $x$-torsion-free left
$R[x,f]$-module $G$. Since $L$ corresponds to the
$G$-special $R$-ideal $\fb$, we must have $L' \subset L$,
since otherwise we would have
$$
(\fb + Rc)R[x,f] \subseteq \grann_{R[x,f]}L' = \grann_{R[x,f]}L =
\fb R[x,f].
$$
Therefore there exists $h \in H$ such that $h + \Gamma_x(H)$ is
annihilated by $\fb R[x,f]$ but not by $(\fb + Rc) R[x,f]$. Since
$\height \fb \geq 1$, it follows from Lemma \ref{tc.1s} that $h
\in 0^*_H$.

Choose a system of parameters $a_1, \ldots, a_d$ for $R$; use
$a_1, \ldots, a_d$ in the notation of \ref{tc.1}(i) for
$H^d_{\fm}(R) = H$, and write $a := a_1\ldots a_d$. There exist $r
\in R$ and $j \in \nn$ such that $h = [r/a^j]$. By
\ref{tc.1r}(vi), we have $r \in ((a_1^j, \ldots, a_d^j)R)^*$.
Since $c$ is a weak parameter test element for $R$, we see that
$cr^{p^n} \in (a_1^{jp^n}, \ldots, a_d^{jp^n})R$ for all $n \gg
0$, so that $cx^nh = 0$ for all $n \gg 0$. Thus there is some $n_0
\in \nn$ such that $cx^n(h + \Gamma_x(H)) = 0$ in $G$ for all $n
\geq n_0$. Therefore, for all $j \in \nn$, we have
$$
x^{n_0}cx^j(h + \Gamma_x(H)) = c^{p^{n_0}}x^{n_0 + j}(h + \Gamma_x(H)) = 0,
$$
so that $cx^j(h + \Gamma_x(H)) = 0$ because
$G$ is $x$-torsion-free. Thus $h + \Gamma_x(H)$ is annihilated by
$RcR[x,f]$ as well as by $\fb R[x,f]$. This is a contradiction.

Therefore $\fb^{[q_0]} \subseteq \sigma'(R) \subseteq \fb$, since $\sigma'(R)$
can be generated by its elements that lie in $R^{\circ}$.
\end{proof}

Use $R'$ to
denote $R$ (as in\/ {\rm \ref{nt.1}}) regarded as an $R$-module by means of $f$.
With this notation, $f : R \lra R'$ becomes a homomorphism
of $R$-modules.
Recall that, when $(R,\fm)$ is a local ring of
dimension $d > 0$, we say that $R$ is
{\em $F$-injective\/} precisely when the induced homomorphisms
$ H^i_{\fm}(f) : H^i_{\fm}(R) \lra H^i_{\fm}(R')$ are injective for all
$i = 0, \ldots, d$. See R. Fedder and K-i. Watanabe \cite[Definition 1.7]{FW87}
and the ensuing discussion.

\begin{cor}
\label{tc.3} Let $(R,\fm)$ be an $F$-injective Cohen--Macaulay
local ring of dimension $d > 0$. The left $R[x,f]$-module $H :=
H^d_{\fm}(R)$ is $x$-torsion-free. By Theorem\/ {\rm \ref{ga.13}},
there exists a (uniquely determined) smallest ideal $\fb$ of
height at least $1$ in the set $\mathcal{I}(H)$ of
$H$-special $R$-ideals.

Let $c$ be any element of $\fb \cap R^{\circ}$. Then $c$ is a
parameter test element for $R$. In fact, $\fb$ is the parameter
test ideal of $R$ (see\/ {\rm \cite[Definition 4.3]{Smith95}}).
\end{cor}

\begin{note} It should be noted that, in Corollary \ref{tc.3}, it is
not assumed that $R$ is excellent.
\end{note}

\begin{proof} With the
notation of Theorem \ref{tc.2}, the HSL-number $m_0$ of $H$ is $0$
when $R$ is $F$-injective, and so $q_0 = 1$ and $G \cong H$ in
this case. By Theorem \ref{tc.2}, each element $c \in \fb \cap
R^{\circ}$ is a $q_0$-weak parameter test element for $R$, that
is, a parameter test element for $R$. Since $R$ has a parameter
test element, its parameter test ideal $\sigma(R)$ is equal to the
ideal of $R$ generated by all parameter test elements. By Theorem
\ref{tc.2}, we therefore have
$
\fb = \fb^{[q_0]} \subseteq \sigma(R) \subseteq \sigma'(R) \subseteq \fb.
$
\end{proof}

\begin{cor}
\label{tc.4} Let $(R,\fm)$ be an $F$-injective Gorenstein local
ring of dimension $d > 0$. The left $R[x,f]$-module $H :=
H^d_{\fm}(R)$ is $x$-torsion-free. By Theorem\/ {\rm \ref{ga.13}},
there exists a (uniquely determined) smallest ideal $\fb$ of
height at least $1$ in the set $\mathcal{I}(H)$ of
$H$-special $R$-ideals.

Let $c$ be any element of $\fb \cap R^{\circ}$. Then $c$ is a test
element for $R$. In fact, $\fb$ is the test ideal of $R$.
\end{cor}

\begin{note} It should be noted that, in Corollary \ref{tc.4}, it is
not assumed that $R$ is excellent.
\end{note}

\begin{proof} This follows immediately from Corollary \ref{tc.3} once
it is recalled that an $F$-injective Cohen--Macaulay local ring is reduced and
that a parameter test element for a reduced Gorenstein local ring $R$
of characteristic $p$ is automatically a test element for $R$: see
the proof of \cite[Proposition 4.1]{Hunek98}.
\end{proof}

\section{\sc Special $R$-ideals and Enescu's $F$-stable primes}
\label{en}

The purpose of this section is to establish connections between
the work in \S \ref{ga} and \S \ref{tc} above and F. Enescu's
$F$-stable primes of an $F$-injective Cohen--Macaulay local ring
$(R,\fm)$, defined in \cite[\S 2]{Enesc03}.

\begin{ntn}
\label{en.1} Throughout this section, $(R,\fm)$ will be assumed to
be a Cohen--Macaulay local ring of dimension $d > 0$, and we shall
let $a_1, \ldots, a_d$ denote a fixed system of parameters for
$R$, and set $\fq := (a_1, \ldots, a_d)R$. We shall use $a_1,
\ldots, a_d$ in the notation of \ref{tc.1}(i) for $H :=
H^d_{\fm}(R)$, and write $a := a_1\ldots a_d$.

For each $b \in R$, we define (following Enescu \cite[Definition
1.1]{Enesc03}) the ideal $\fq(b)$ by
$$
\fq(b) := \left\{ c \in R : cb^{p^n} \in \fq^{[p^n]} \mbox{~for
all~} n \gg 0 \right\}.
$$
(Actually, Enescu only made this definition when $b \not\in \fq$;
however, the right-hand side of the above display is equal to $R$
when $b \in \fq$, and there is no harm in our defining $\fq(b)$ to
be $R$ in this case.) In view of \ref{tc.1}(i), the ideal $\fq(b)$
is equal to the ultimate constant value of the ascending chain
$(\fb_n)_{n \in \nn}$ of ideals of $R$ for which
$\bigoplus_{n\in\nn}\fb_n x^n = \grann_{R[x,f]}R[x,f][b/a]$, the
graded annihilator of the $R[x,f]$-submodule of $H$ generated by
$[b/a]$.

Now consider the special case in which $R$ is (also)
$F$-injective. Then the left $R[x,f]$-module $H$ is
$x$-torsion-free, and so it follows from Lemma \ref{nt.5} that,
for each $b \in R$, the ideal $\fq(b)$ is radical and
$\grann_{R[x,f]}R[x,f][b/a] = \fq(b)R[x,f]$; thus $\fq(b)$ is an
$H$-special $R$-ideal. We again follow Enescu and set
$$
Z_{\fq,R} := \{ \fq(b) : b \in R \setminus \fq \}.
$$
\end{ntn}

Enescu proved, in \cite[Theorem 2.1]{Enesc03}, that (when
$(R,\fm)$ is Cohen--Macaulay and $F$-injective) the set of maximal
members of $Z_{\fq,R}$ is independent of the choice of $\fq$, is
finite, and consists of prime ideals. The next theorem shows that
the set of maximal members of $Z_{\fq,R}$ is actually equal to the
set of maximal members of $\mathcal{I}(H) \setminus \{R\}$: we saw
in Lemma \ref{ga.11} that this set consists of prime ideals, and
in Corollary \ref{ga.12c} that it is finite.

\begin{thm}
\label{en.2} Let the situation and notation be as in\/ {\rm
\ref{en.1}}, and suppose that the Cohen--Macaulay local ring
$(R,\fm)$ is $F$-injective. Then the set of maximal members of
$Z_{\fq,R}$ is equal to the set of maximal members of
$\mathcal{I}(H) \setminus \{R\}$.
\end{thm}

\begin{proof} The comments in \ref{en.1}
show that $Z_{\fq,R} \subseteq \mathcal{I}(H)$; clearly, no member
of $Z_{\fq,R}$ can be equal to $R$. It is therefore sufficient for
us to show that a maximal member $\fp$ of $\mathcal{I}(H)
\setminus \{R\}$ must belong to $Z_{\fq,R}$.

Let $L \in \mathcal{A}(H)$ be the special annihilator submodule of
$H$ corresponding to $\fp$. Now $H$ is an Artinian $R$-module: let
$h$ be a non-zero element of the socle of $L$. By Lemma
\ref{ga.11}, we have $\grann_{R[x,f]}R[x,f]h = \fp R[x,f]$.
However, for each $j \in \N$, we have $R[1/a^j] \cong R/(a_1^j,
\ldots, a_d^j)R$, by \ref{tc.1r}(vii), so that
$$
\Hom_R(R/\fm,R[1/a^j]) \cong \Hom_R(R/\fm,R/(a_1^j, \ldots, a_d^j)R) \cong
\Ext^d_R(R/\fm,R)
$$
by \cite[Lemma 1.2.4]{BH}. It follows that it is possible to write $h$ in
the form $h = [r/a]$ for some $r \in R$, and therefore $\fp =
\fq(r) \in Z_{\fq,R}$.
\end{proof}

\bibliographystyle{amsplain}

\end{document}